\documentclass[12pt, reqno, a4paper,oneside]{amsart}
\usepackage{graphicx, epsfig, psfrag}
\usepackage{amsfonts,amsmath,amssymb,amsbsy,amsthm}
\usepackage{bm}
\usepackage{color}
\usepackage{float}
\usepackage{mathrsfs}
\usepackage{longtable}
\usepackage[top=3cm,bottom=3cm,left=3cm,right=3cm]{geometry}
\usepackage{epstopdf}
\usepackage[normalem]{ulem}

\usepackage{environments}
\numberwithin{theorem}{section}
\numberwithin{equation}{section}

\renewcommand{\paragraph}[1]{\subsubsection{#1}}

\renewcommand{\cases}[1]{\left\{ \begin{array}{rl} #1 \end{array} \right.}

\def\XXint#1#2#3{{\setbox0=\hbox{$#1{#2#3}{\int}$ }
\vcenter{\hbox{$#2#3$ }}\kern-.6\wd0}}

\def\b{\big}
\def\B{\Big}
\def\bg{\bigg}

\def\sep{\,|\,}
\def\bsep{\,\b|\,}

\def\R{\mathbb{R}}

\def\Z{\mathbb{Z}}

\def\dx{\,{\rm d}x}

\def\pp{\partial}

\def\<{\langle}
\def\>{\rangle}

\def\mA{{\sf A}}
\def\mB{{\sf B}}

\def\mF{{\sf F}}

\def\mQ{{\sf Q}}

\def\D{\nabla}
\def\del{\delta}
\def\ddel{\delta^2}

\def\a{{\rm a}}
\def\c{{\rm c}}
\def\ac{{\rm ac}}
\def\i{{\rm i}}

\def\L{\Lambda}

\def\Us{\mathscr{W}}

\def\Ush{\dot{\Us}^{1,2}}

\def\yF{y_\mF}

\def\E{E}
\def\Ea{\E^\a}

\def\Ei{\E^\i}
\def\Ec{\E^\c}

\def\Eac{\E^\ac}

\def\L{\Lambda}
\def\La{\L^{\a}}
\def\Li{\L^{\i}}
\def\Lc{\L^{\c}}
\def\yd{y_0}
\def\Nhd{\mathscr{N}}
\def\rcut{r_{\rm cut}}
\def\Rg{\mathscr{R}}
\def\Rgp{\Rg^{+}}
\def\vsig{\varsigma}
\def\T{\mathscr{T}_h}
\def\UsT{\Us_h}

\def\Tmu{\mathscr{T}_\mu}
\def\vor{\rm vor}

\def\defc{\rm def}
\def\shift{\rm shift}

\definecolor{cocol}{rgb}{0.7, 0, 0}
\definecolor{lzcol}{rgb}{0, 0, 0.9}
\definecolor{todocol}{rgb}{0.0, 0.4, 0.0}

\begin{document}

\title[Energy-Based A/C Coupling Without Ghost Forces]{Energy-Based
  Atomistic-to-Continuum \\ Coupling Without Ghost Forces}

\author{C. Ortner}
\address{C. Ortner\\ Mathematics Institute \\ Zeeman Building \\
  University of Warwick \\ Coventry CV4 7AL \\ UK}
\email{christoph.ortner@warwick.ac.uk}

\author{L. Zhang}
\address{L. Zhang \\ Department of Mathematics,
  Institute of Natural Sciences, and MOE Key Lab in Scientific and Engineering Computing\\
  Shanghai Jiao Tong University \\ 800 Dongchuan Road \\ Shanghai
  200240 \\ China}
\email{lzhang@sjtu.edu.cn}

\date{\today}

\thanks{CO's work was supported by EPSRC grant EP/J022055/1, and by
  the Leverhulme Trust. LZ's work was supported by One Thousand Plan
  of China for young scientists. Parts of this work were developed
  during the 2013 IPAM semester programme `Materials Defects:
  Mathematics, Computation, and Engineering'.}


\begin{abstract}
  We present a practical implementation of an energy-based
  atomistic-to-continuum (a/c) coupling scheme without ghost forces,
  and numerical tests evaluating its accuracy relative to other types
  of a/c coupling schemes.
\end{abstract}

\maketitle



\section{Introduction}
\label{sec:intro}
Atomistic-to-continuum coupling methods (a/c methods) are a class of
computational multiscale schemes that combine the accuracy of
atomistic models of defects with the computational efficiency of
continuum models of elastic far-fields \cite{Gumbsch:1989,
  Ortiz:1995a,XiBe:2004,Shimokawa:2004}. In the present article, we
present the first succesful implementation of a practical {\em patch
  test consistent} energy based a/c coupling scheme. Previously such
schemes were only available for 2-body interactions
\cite{Shapeev:2010a, Shapeev2012}

In recent years a numerical analysis theory of a/c methods has
emerged; we refer to \cite{LuOr:acta} for a review.  This theory has
identified three prototypical classes of a/c schemes: patch test
consistent energy-based coupling, force-based coupling (including
force-based blending), and energy-based blending. The classical
numerical analysis concepts of consistency and stability are applied
to precisely quantify the errors committed in these schemes, and clear
guidelines are established for their practical implementation
including optimisation of approximation parameters. The results in
\cite{Dobson:2008b,emingyang,LuOr:acta,Or:2011a,PRE-ac.2dcorners}
indicate that patch test consistent a/c couplings observe
(quasi-)optimal error estimates in the energy-norm. However, to this
date, no general construction and implementation of such schemes has
been presented. Instead, one normally compromises by either turning to
patch test consistent force-based schemes
\cite{Shenoy:1999a,MiLu:2011,2013-bqcfcomp, Gumbsch:1989} or to
blending schemes \cite{XiBe:2004,2012-optbqce} which have some control
over the consistency error. Quasi-optimal implementations of such
schemes are described in \cite{2012-optbqce, 2013-bqcfcomp}.

Existing patch test consistent schemes are restricted in their range
of validity: \cite{Shimokawa:2004} is only consistent for flat a/c
interfaces and short-ranged interactions, \cite{E:2006} extends the
idea to arbitrary range and \cite{PRE-ac.2dcorners} to domains with
corners (but restricting again to nearest-neighbour interaction). On
the other hand, the schemes presented in \cite{Shapeev:2010a,
  Shapeev2012, Makridakis2012} are valid for general interaction range
and a/c interfaces with corners, but are restricted to pair
interactions.

In the present article, we shall present a generalisation of the
geometric reconstruction technique \cite{Shimokawa:2004, E:2006,
  PRE-ac.2dcorners}, which we subsequently denote GRAC. Briefly, the
idea is that, instead of evaluating the interatomic potential near the
a/c interface with atom positions obtained by interpolating the
continuum description, one extrapolates atom positions from those in
the atomistic region (geometric reconstruction). This idea is somewhat
analogous to the implementation of Neumann boundary conditions for
finite difference schemes. There is substantial freedom in how this
reconstruction is achieved, leading to a number of {\em free
  parameters}. One then determines these {\em reconstruction
  parameters} by solving the ``geometric consistency equations''
\cite{E:2006}, which encode a form of patch test consistency and lead
to a first-order consistent coupling scheme \cite{Or:2011a}.

The works \cite{E:2006, PRE-ac.2dcorners, Or:2011a} have demonstrated
that GRAC is a promising approach, but also indicate that explicit
analytical determination of the reconstruction parameters for general
a/c interface geometries with general interaction range may be
impractical. Instead we propose to compute the reconstruction
parameters in a preprocessing step. Although this is a natural idea it
has not been pursued to the best of our knowledge. 

A number of challenges must be overcome to obtain a robust numerical
scheme in this way. The two key issues we will discuss are:
\begin{itemize}
\item[(A)] If the geometric consistency equations have a solution then
  it is not unique. The consistency analysis \cite{Or:2011a} suggests
  that a solution is best selected through $\ell^1$-minimisation of
  the coefficients. Indeed, we shall demonstrate that the least
  squares solution leads to prohibitively large errors.
\item[(B)] In \cite{2013-stab.ac} we proved that there exists no
  universally stable a/c coupling of geometric reconstruction type. We
  will see that this is in fact of practical concern and demonstrate
  that the stabilisation mechanism proposed in \cite{2013-stab.ac}
  appears to resolve this issue.
\end{itemize}

In the remainder of the paper we present a complete description of a
practical implementation of the GRAC method (\S~\ref{sec:formulation})
and numerical experiments focused primarily on investigating
approximation errors (\S~\ref{sec:numerical_tests}). We will comment
on open issues and possible improvements in \S~\ref{sec:discuss},
which are primarily concerned with the computational cost of
determining the reconstruction coefficients.

\section{Formulation of the GRAC Method}
\label{sec:formulation}
In formulating the GRAC scheme, we adopt the point of view of
\cite{2013-defects}, where the computational domain and boundary
conditions are considered part of the approximation. This setting is
convenient to assess approximation errors. Adaptions of the coupling
mechanism to other problems are straightforward.

We first present a brief review, ignoring some technical details, of a
model for crystal defects in an infinite lattice from
\cite{2013-defects}, and some results concerning their structure
(\S~\ref{sec:formulation:atm}). In \S~\ref{sec:formulation:ac} we
present a generic form of a/c coupling schemes, which we then
specialize to the GRAC scheme in \S~\ref{sec:formulation:grac}. In
\S~\ref{sec:formulation:grac} and in \S~\ref{sec:optim_coeffs} we
address, respectively, the two key issues (A) and (B) mentioned in the
introduction.

For the sake of simplicity of presentation, and to emphasize the
algorithmic aspects of the GRAC method, we restrict the presentation
to relatively simple settings such as point defects and microcracks as
in \cite{2012-optbqce, 2013-bqcfcomp}. The concepts required to
generalize the presentation to problems involving dislocations can be
found in \cite{2013-defects}.

\subsection{Atomistic model}
\label{sec:formulation:atm}
\def\Omdef{\Omega^{\rm def}}
Let $d \in \{2, 3\}$ denote the problem dimension. Fix a non-singular
$\mA \in \R^{d \times d}$ to define a Bravais lattice $\mA \Z^d$. Let
$\L \subset \R^d$ be a discrete reference configuration of a crystal,
possibly with a local defect: for some compact domain $\Omdef$ we
assume that $\L \setminus \Omdef = \mA\Z^d \setminus \Omdef$ and $\L
\cap \Omdef$ is finite. It can be readily seen \cite{2013-defects},
that certain point defects (e.g., interstials, vacancies; see Figure
\ref{fig:point_defects}) can be enforced that way.

\definecolor{grey}{rgb}{.5, .5, .5}
\begin{figure}
  \begin{minipage}{4.5cm}
    \includegraphics[height=4cm]{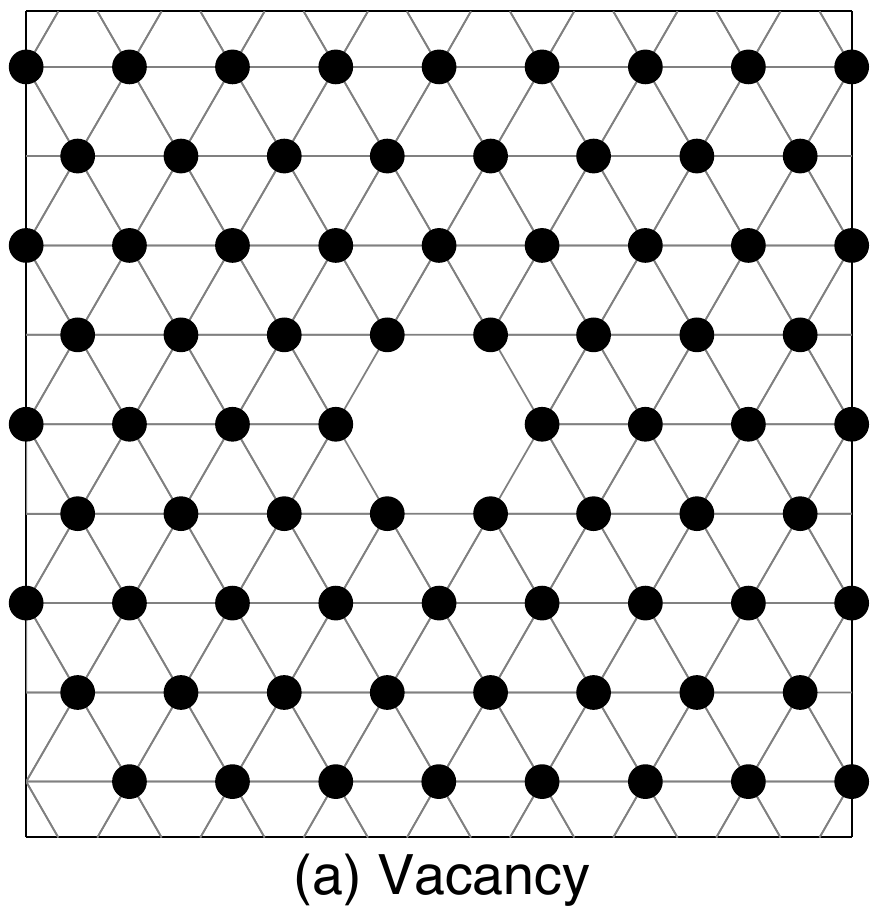} \\[2mm]
    \includegraphics[height=4cm]{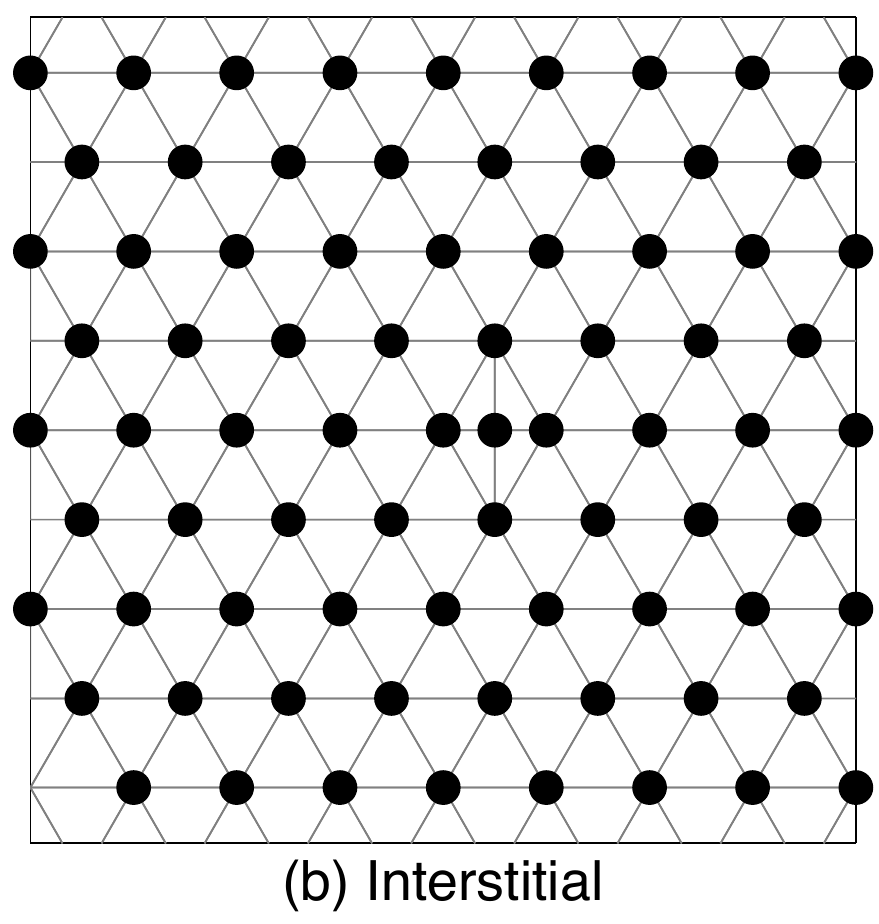}
  \end{minipage}
  \begin{minipage}{9cm}
    \includegraphics[height=7.5cm]{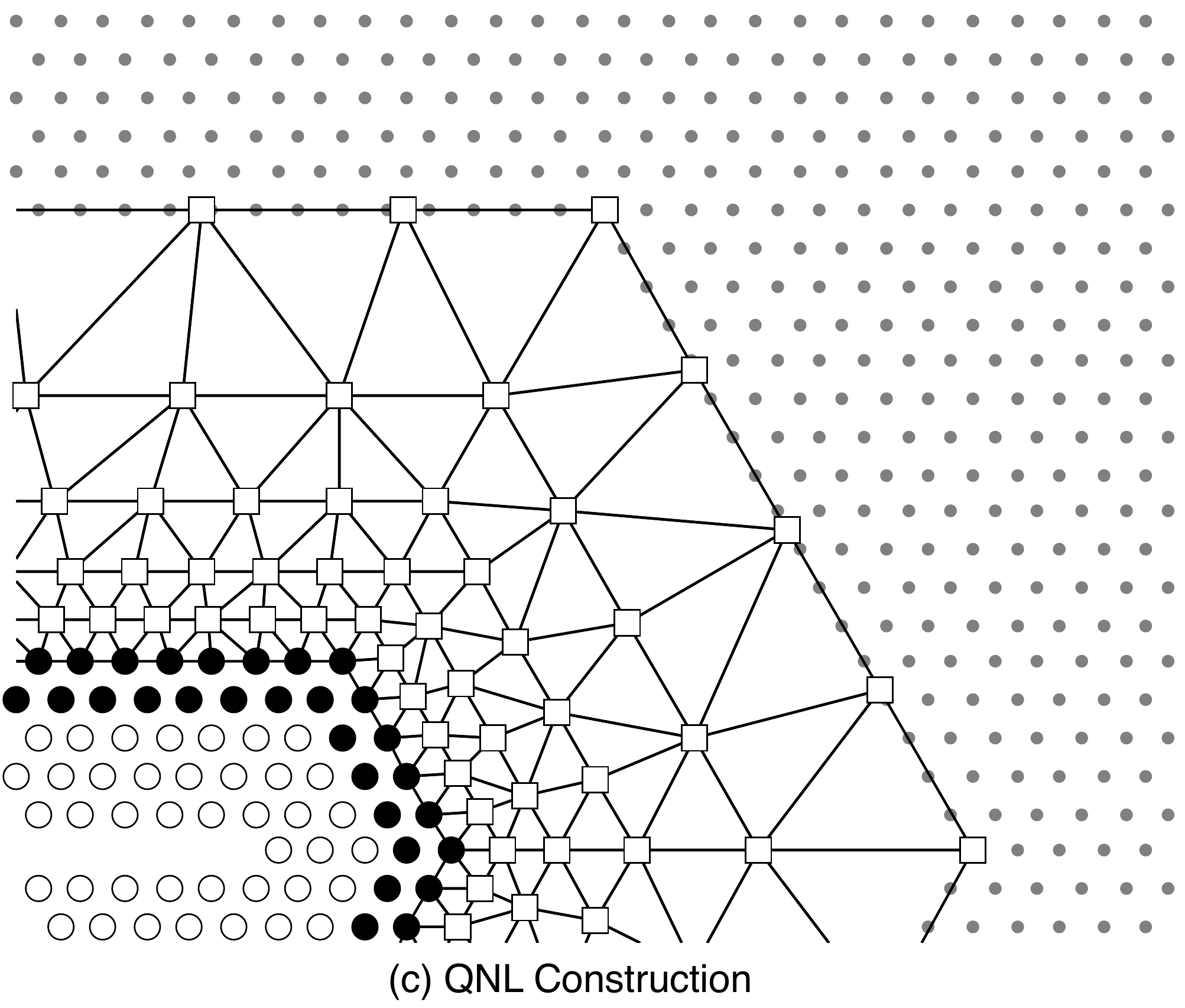}
  \end{minipage}
  \caption{\label{fig:point_defects} (a, b) Examples of reference
    configurations $\L$ with point defects embedded. (c) Construction
    of the QNL method: $\circ$ atomistic potential $\Phi_\ell$;
    $\bullet$ interface potential $\Phi_\ell^\i$; $\square$
    Cauchy--Born potential $W$ (precisely, $W$ is applied on
    elements/triangles); {\color{grey} \footnotesize $\bullet$}
    far-field boundary condition $y(\ell) = \mB \ell$ is imposed.}
\end{figure}


To avoid minor technical difficulties, we prescribe a maximal
interaction neighbourhood in the reference configuration. This is a
restriction that can be lifted with little additional work
\cite[Remark 2.1]{2013-defects}. For each $\ell \in \L$ we denote this
neighbourhood by $\Nhd(\ell) := \{ \ell' \in \L \sep |\ell'-\ell| \leq
\rcut \}$, for some specified reference cut-off radius $\rcut$. We
define the assocated sets $\Nhd_*(\ell) := \Nhd(\ell) \setminus
\{\ell\}$ and $\Rg(\ell) := \{ \ell'-\ell \sep \ell' \in
\Nhd_*(\ell)\}$. We define the ``finite difference stencil'' $Dv(\ell)
:= (D_\rho v(\ell))_{\rho \in \Rg(\ell)} :=
(v(\ell+\rho)-v(\ell))_{\rho \in \Rg(\ell)}$. Higher-order finite
differences, $D_\rho D_\vsig v$ and $D^2 v$ are defined in a canonical
way.

We use this notation to define a discrete energy space. For $v : \L
\to \R^m$, let the discrete energy-norm be defined by
\begin{displaymath}
  \| v \|_{\Ush} := \| Dv \|_{\ell^2} := \bg( \sum_{\ell \in \L}
  \sum_{\rho \in \Rg(\ell)} \frac{|D_\rho v(\ell)|^2}{|\rho|^2}
  \bg)^{1/2}
  = \bg( \sum_{\ell \in \L} \sum_{\ell' \in \Nhd_*(\ell)} 
  \frac{|v(\ell')-v(\ell)|^2}{|\ell'-\ell|^2} \bg)^{1/2},
\end{displaymath}
which we can think of as a discrete $H^1$-seminorm.
%
%
%
Then, the associated discrete function space is defined by
\begin{align*}
  \Ush &:= \b\{ u : \L \to \R^m \bsep \| u \|_{\Ush} < +\infty \b\}.
\end{align*}
The space $\Ush$ can be thought of as the space of all relative
displacements with finite energy.

For a deformed configuration $y : \L \to \R^d$ and $\ell \in \L$, let
$\Phi_\ell(y) = \Phi_\ell((y_{\ell'})_{\ell' \in \Nhd(\ell)})$ denote
a site energy functional associated with $\ell$. For $\ell \in \L
\setminus \Omdef$ we assume that $\Phi_\ell(y) \equiv \Phi(y -
y(\ell))$, i.e., the crystal is homogeneous outside $\Omdef$. By
changing the interaction potential inside $\Omdef$, impurities or
``cut bonds'' can be modelled.

The prototypical example is the embedded atom model \cite{Daw:1984a},
for which $\Phi_\ell$ is of the form 
\begin{equation}
  \label{eq:eam}
  \begin{split}
  \Phi_\ell(y) &= \sum_{\ell' \in \Nhd_*(\ell)} \phi\b(|y(\ell')-y(\ell)|\b) +
  F\B( {\textstyle \sum_{\ell' \in \Nhd_*(\ell)}} \psi\b( |y(\ell') -
  y(\ell)|\b) \B) \\
  &= \sum_{\rho \in \Rg(\ell)} \phi\b(|D_\rho y(\ell)|\b) + F\B(
  {\textstyle \sum_{\rho \in \Rg(\ell)}} \psi\b( |D_\rho y(\ell)|\b)
  \B). 
\end{split}
\end{equation}
The energy of an infinite configuration is typically ill-defined, but
the {\em energy-difference functional}
\begin{displaymath}
  E(y; z) = \sum_{\ell \in \L} \Phi_\ell(y) - \Phi_\ell(z)
\end{displaymath}
\def\yB{y^\mB}
is a meaningful object. For example, if $y - z$ has compact support,
then $E(y;z)$ is well-defined. More generally it is shown in
\cite[Thm. 2.2]{2013-defects}, under natural technical conditions on
the site potentials $\Phi_\ell$, that $u \mapsto E(\yB+u; \yB), u \in
\Ush$, is well-defined and (Fr\'echet) differentiable, where $\yB(x) =
\mB x$.


Given a {\em macroscopic applied strain} $\mB \in \R^{d \times d}$, we
aim to compute
\begin{equation}
  \label{eq:min}
  y \in \arg\min \b\{ E(y; \yB) \bsep y - \yB \in \Ush \b\}.
\end{equation}
A solution to \eqref{eq:min} will satisfy the far-field boundary
condition $y(\ell) \sim \mB \ell$ as $|\ell| \to \infty$, imposed
through the condition that $y - \yB \in \Ush$.

We call a solution $y$ {\em strongly stable} if there exists $c_0 > 0$
such that $\< \ddel E(y) v, v \> \geq c_0 \| Dv \|_{\ell^2}^2$ for all
$v \in \Ush$.

Here, and throughout, we write $\del^j E(y)$ instead of $\del^j E(y;
z)$ since the variations of the energy difference only depend on the
first component.

\begin{remark}
  The far-field boundary condition $y(\ell) \sim \mB \ell$ can be
  generalised to any deformation $\yd$ satisfying $\del E(\yd) \in
  (\Ush)^*$, for example, to dislocations by replacing $\mB \ell$ with
  the linear elasticity solution of the dislocation
  \cite{2013-defects}.
  %
\end{remark}

\subsection{A/C coupling}
\label{sec:formulation:ac}
We begin by giving a generic formulation of an a/c coupling, which we
subsequently make concrete employing concepts and notation from
various earlier works, such as
\cite{Ortiz:1995a,Shenoy:1999a,Shimokawa:2004,2012-optbqce}, but
adapting the formulation to our setting of \S~\ref{sec:formulation}.
The construction is visualised in Figure \ref{fig:point_defects}(c).


To choose a computational domain let $\Omega \subset \R^d$ be a simply
connected, polygonal and closed set. We decompose $\Omega = \Omega^\a
\cup \Omega^\c$, where $\Omega^\a$ is again simply connected and
polygonal, and contains the defect: $\Omdef \subset \Omega^\a$. Let
$\T$ be a regular partition of $\Omega^\c$ into triangles ($d = 2$) or
tetrahedra ($d = 3$). Let $I_h$ denote the associated nodal
interpolation operator.

Next, we decompose the set $\L^{\a,\i} := \L \cap \Omega^\a = \L^\a
\cup \L^\i$ into a core atomistic region $\L^\a$ and an interface
region $\L^\i$ (typically a few ``layers'' of atoms surrounding
$\L^\a$).

We can now define the space of {\em coarse-grained} displacement maps,
\begin{align*}
  \UsT := \b\{ u_h : \Omega^\c \cup \L^{\a,\i} \to \R^m \bsep ~&
  \text{ $u_h$ is continuous and p.w. affine w.r.t. $\T$, } \\[-1mm]
  & \text{ and $u_h = 0$ on $\partial \Omega$ } \b\}.
\end{align*}
The associated space of coarse-grained deformations is $\yB + \UsT$.

The Cauchy--Born strain energy function is given by
\begin{displaymath}
  W(\mF) := |{\rm vor}(\ell)|^{-1} \Phi_\ell( \mF \cdot \Rg(\ell) ) \qquad \text{ for some
  } \ell  \in \L \setminus \Omdef,
\end{displaymath}
where ${\rm vor}(\ell)$ is the voronoi cell associated with $\ell$.
(Due to the homogeneity of the lattice and interaction outside
$\Omdef$, the definition is independent of $\ell$.)

For $\ell\in\L^{\i}$, we choose a modified interface site potential
$\Phi_\ell^\i$ and an effective cell $v_\ell^\i \subset {\rm
  vor}(\ell)$ associated with $\ell$ (specific choices will be
specified in \S~\ref{sec:formulation:grac}), and define the effective
volume associated with $\ell$ as
$\omega^{\i}_\ell:= |v_\ell^\i| / |{\rm vor} (\ell)|$. Further, for
each element $T \in \T$ we define the effective volume
  $\omega_T := |T \setminus (\cup_{\ell \in \L^{\i}} v^\i_\ell)|$.

Then, a generic a/c
coupling energy difference functional is then defined by
\begin{align}
  \label{eq:generic_ac_energy}
  E^\ac(y_h; z_h) &:= \sum_{ \ell \in \L^\a} \B( \Phi_\ell(y_h) -
  \Phi_\ell( z_h) \B) + \sum_{\ell \in \L^\i} \omega^\i_\ell\B( \Phi^\i_\ell(y_h) -
    \Phi_\ell^\i(z_h) \B) \\
    \notag
    & \qquad + \sum_{T \in \T} \omega_T \B( W(\D y_h|_T) -
    W(\D z_h|_T) \B). 
\end{align} 
Thus, we obtain the approximate variational problem
\begin{equation}
  \label{eq:min_ac}
  y_h \in \arg\min \b\{ E^\ac(y_h; \yB) \bsep y_h - \yB \in
  \UsT \b\}.
\end{equation}

\subsubsection{The patch tests}
\def\yF{y^\mF}
A key condition that has been widely discussed in the a/c coupling
literature is that $E^\ac$ should exhibit no ``ghost
forces''. Following the language of \cite{Or:2011a}, we call this
condition the {\em force patch test}: for $\L = \mA\Z^d$ and
$\Phi_\ell \equiv \Phi$ (homogeneous
lattice without defects)
\begin{equation}
  \label{eq:force_pt}
  \< \del E^\ac(\yF), v \> = 0 \qquad \forall v \in \UsT,
  \quad \mF \in \R^{d \times d}.
\end{equation}
In addition, to guarantee that $E^\ac$ approximates the atomistic
energy $E$, it is reasonable to also require that the interface
potentials satisfy an {\em energy patch test}
\begin{equation}
  \label{eq:energy_pt}
  \Phi_\ell^\i(\yF) = \Phi_\ell(\yF) \qquad \forall \mF \in \R^{d
    \times d}, \quad \ell \in \L^\i.
\end{equation}


\subsection{General GRAC formulation}
\label{sec:formulation:grac}
To complete the definition of the a/c coupling energy
\eqref{eq:generic_ac_energy} and of the associated variational problem
\eqref{eq:min_ac}, we must specify the interface region $\L^\i$, the
interface site potentials $\Phi_\ell^\i$ and the associated volumes
$\omega_\ell^\i$. The approach we present here is an extension of
\cite{Shimokawa:2004, E:2006, PRE-ac.2dcorners}.

First we note that, due to homogeneity of $\Phi_\ell$ outside of
$\Omdef$, we can write 
\begin{displaymath}
  \Phi_\ell(y) = V\b( D y(\ell) \b),
\end{displaymath}
for some potential $V$ that is a function of the finite differences
instead of a function of positions.

We now define $\Phi_\ell^\i$ in terms of $V$.  For each $\ell \in
\L^\i, \rho, \vsig \in \Rg(\ell)$, we let $C_{\ell;\rho,\vsig}$ be
free parameters, and define
\begin{equation}
  \label{eq:defn_Phi_int}
  \Phi_\ell^\i(y_h) := V \B( \b( {\textstyle \sum_{\vsig \in
      \Rg(\ell)} C_{\ell;\rho,\vsig} D_\vsig y_h(\ell) } \b)_{\rho \in
    \Rg(\ell)} \B)
\end{equation}
A convenient short-hand is
\begin{displaymath}
  \Phi_\ell^\i(y_h) = V( C_\ell \cdot Dy_h(\ell) ) \quad \text{where}
  \quad \cases{
    C_\ell := (C_{\ell;\rho,\vsig})_{\rho,\vsig \in \Rg(\ell)}, \quad
    \text{and} &\\
    C_\ell \cdot Dy := \b( {\textstyle \sum_{\vsig \in
      \Rg(\ell)} C_{\ell;\rho,\vsig} D_\vsig y } \b)_{\rho \in
    \Rg(\ell)}. & }
\end{displaymath}
We call $C_{\ell;\rho,\vsig}$ the {\em reconstruction
  parameters}.

The parameters are to be chosen so that the resulting energy
functional $E^\ac$ satisfies the energy and force patch tests
\eqref{eq:force_pt} and \eqref{eq:energy_pt}. 

\begin{remark}
  \label{rem:terminology_grac}
  The approach \eqref{eq:defn_Phi_int} is labelled {\em quasi-nonlocal
    coupling} in \cite{Shimokawa:2004} since the coefficients are
  (typically) chosen so that the interaction of $\L^\i$ with the
  atomistic region $\L^\a$ is non-local while the interaction of
  $\L^\i$ with the continuum region is local. In \cite{E:2006} the
  approach is labelled {\em geometric reconstruction} since we can
  think of the operation $C_\ell \cdot Dy(\ell)$ as {\em
    reconstructing} atom positions in the continuum region, using only
  information from the atomistic region and interface. 

  A more pragmatic point of view is to simply view the atomistic model
  and continuum model as two different finite difference schemes for
  the same PDE and to ``fit'' parameters that would {\em consistently
    patch them together}.
  %
\end{remark}

\subsubsection{Energy patch test} 
A sufficient and necessary condition for the energy patch test
\eqref{eq:energy_pt} is that $\mF \cdot \Rg(\ell) = C_\ell \cdot (\mF
\cdot \Rg)$ for all $\mF \in \R^{m \times d}$ and $\ell \in
\L^\i$. This is equivalent to
\begin{equation}
  \label{eq:E_pt_grac}
  \rho = \sum_{\vsig \in \Rg(\ell)} C_{\ell; \rho,\vsig} \vsig \qquad
  \forall \ell \in \L^\i, \quad \rho \in \Rg(\ell).
\end{equation}


\subsubsection{Force patch test}
\label{sec:force_pt_lineqn}
The force patch test \eqref{eq:force_pt} leads to a fairly complex set
of equations. From the general GRAC formulation
\eqref{eq:generic_ac_energy}, we can decompose the first variation of
the A/C coupling energy into three parts,
\begin{displaymath}
  \<\del\Eac(\yF), u\> = \<\del\Ea(\yF), u\> 
  + \<\del\Ei(\yF), u\> + \<\del\Ec(\yF), u\>.
\end{displaymath}
To simplify the notation, we drop the $\yF$ dependence from the
expression, for example, we write $\Ea$ instead of $\Ea(\yF)$,
$\D_\rho V$ instead of $\D_\rho V(D\yF)$, and so forth. Here, $\D_\rho
V$ denotes the partial derivative of $V$ with respect to the $D_\rho y$
component.

Since $\D_\rho V = -\D_{-\rho} V$, we only consider half of the
interaction range: we fix $\Rgp \subset \Rg$ such that $\Rgp \cup
(-\Rgp) = \Rg$ and $\Rgp \cap (-\Rgp) = \emptyset$.

The first variations in the a/c coupling energy can be expanded into
the following expressions,
\begin{align*}
  \<\del\Ea, u\> &= \sum_{\substack{\rho \in \Rgp \\ \ell \in
      \La-\rho}} {\b[\D_\rho V \cdot u(\ell) \b]}- \sum_{\substack{\rho \in
      \Rgp \\ \ell\in\La+\rho}}\b[ \D_\rho V \cdot u(\ell)\b], \\
  \<\del\Ei, u\>  &= \sum_{\substack{\vsig \in \Rg \\ \ell \in
      \Li+\vsig}} \omega_{\ell-\vsig}^\i \sum_{\rho\in\Rgp}
  (C_{\ell-\vsig;\rho,\vsig}-C_{\ell-\vsig;-\rho,\vsig}) \b[\D_\rho V \cdot u(\ell)\b]\\
  & \qquad
  -\sum_{\ell\in\Li} \omega_\ell^\i \sum_{\rho\in\Rgp}
  \sum_{\vsig\in\Rg} (C_{\ell;\rho,\vsig}-C_{\ell;-\rho,\vsig}) \b[\D_\rho
  V \cdot u(\ell) \b], \quad \text{and} \\
  \<\del\Ec, u\> &= \sum_T\sum_{\rho\in\Rgp}\sum_{i=1}^3
  2\frac{\omega_T}{|\vor|}\D_T\phi_i^T\cdot \rho \b[ {\D_\rho V} \cdot
  u^T_i \b],
\end{align*}
where the nodes $\ell^T_i$ are the three corners of the triangle $T$,
$u^T_i = u(\ell^T_i)$ and $\phi_i^T$ are the three nodal linear basis
corresponding to $u^T_i$, $i=1,2,3$. The complete calculations are
shown in \S~\ref{sec:app_delEac}.

Since we require that the force patch test \eqref{eq:force_pt} holds
for {\em all potentials} $V$, we can think of $\D_\rho V \cdot
u(\ell)$ as independent symbols.  Collecting all the coefficients for
the terms $\D_\rho V \cdot u(\ell)$, we obtain
\begin{align*}
	\<\del \Ea, u\> & = \sum_{\ell\in\La +\Rg}\sum_{\rho\in\Rgp}
        c^\a_\rho(\ell) \b[\D_\rho V \cdot u(\ell)\b] \\
	\<\del \Ei, u\> & = \sum_{\ell\in\Li
          +\Rg}\sum_{\rho\in\Rgp}c^\i_\rho(\ell) \b[ \D_\rho V \cdot
        u(\ell) \b] \\
	\<\del \Ec, u\> & =
        \sum_{\ell\in\Lc}\sum_{\rho\in\Rgp}c^\c_\rho(\ell) \b[ \D_\rho
        V \cdot u(\ell)   \b].
\end{align*}
%
The coefficients $c^\a_\rho(\ell)$, {$c^\i_\rho(\ell)$} and
$c^\c_\rho(\ell)$ are geometric parameters of the underlying lattice
and of the interface geometry, while the coefficients
$c^\i_\rho(\ell)$ also dependend linearly on the unknown
reconstruction paramters $C_{\ell;\rho,\vsig}$.

Since force patch test is automatically satisfied for the atomistic
model and the Cauchy--Born continuum model, we only need to consider
the force consistency for those sites which the modified interfacial
potential can influence, namely, the extended interface region
$\L^\i+\Rg:=\{\ell\in\L|\exists \ell'\in\L^\i, \exists\rho\in\Rg,\text{ such that } \ell = \ell'+\rho\}$. 

We summarize the foregoing calculation in the following result.

\begin{proposition}
  \label{prop:coefficients}
  A necessary and sufficient condition on the reconstruction
  parameters $C_\ell$ to satisfy the force patch test
  \eqref{eq:force_pt}  for all $V \in C^\infty((\R^d)^\Rg)$ is
  \begin{equation}
    \label{eq:F_pt_grac}
    c^\a_\rho(\ell) + c^\i_\rho(\ell) + c^\c_\rho(\ell) = 0
  \end{equation}
  for $\ell\in\Li+\Rg$, and $\rho\in\Rgp$.
  %
\end{proposition}

\bigskip

\begin{figure}
  \begin{center}
    \includegraphics[width=6cm]{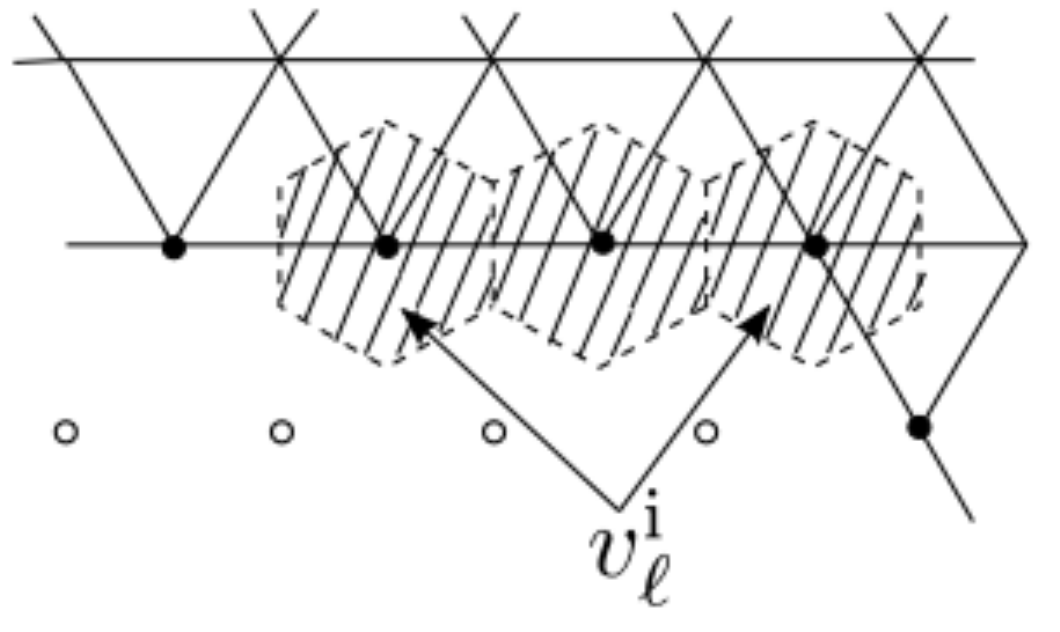}
    \includegraphics[width=6cm]{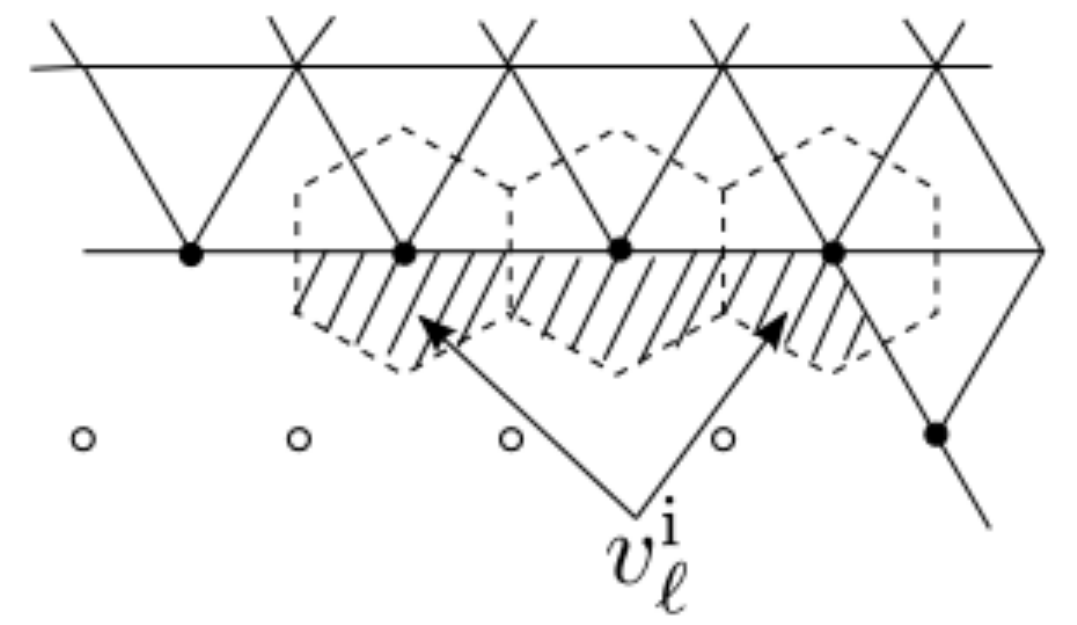}
    \caption{\label{fig:eff_vol} Effective Voronoi cells for the
      interface nodes (filled circles) are the shaded area in the
      above figure. Left figure corresponds to METHOD 1, right figure
      corresponds to METHOD 2. Different choices of effective cells
      result in different values of weights $\omega^{\i}_{\ell}$, for
      method 1, $\omega_\ell^\i = 1$, and for method 2,
      $\omega_\ell^\i<1$ for the outmost interface atoms which are
      adjacent to the continuum region.}\label{fig:effvol}
  \end{center}
\end{figure}

At this stage there is still some freedom in the design of GRAC type
a/c couplings. We implemented the following two variants which place
some additional restrictions, but still do not fully define the
method. See also Figure \ref{fig:eff_vol}.

\begin{itemize}
\item {\bf METHOD 1} is an extension of the construction in
  \cite{PRE-ac.2dcorners}. We choose $v_\ell^\i = {\rm vor}(\ell)$ for
  all $\ell \in \Li$. No other constraints are placed on the method.

  For practical purposes, this method normally requires that in
  $\Omega^\c$, within several layers of atoms surrounding $\L^\i$ all
  nodes of the finite element mesh precisely coincide with the atomic
  sites in these layers; see \S~\ref{sec:appendix:grac}.

\item {\bf METHOD 2} is a variation and extension of the local
  reflection method that is briefly discussed in
  \cite{2013-stab.ac}. We choose $v_\ell^\i = {\rm vor}(\ell) \cap
  \Omega^\a$, and also constrain
  \begin{displaymath}
    C_{\ell; \rho,\vsig} = 0 \quad \text{ for } {\ell \in \Li, \ell+\vsig\in \Omega^\c}.
  \end{displaymath}
  This has the advantage that we now only need to impose the force
  balance equation for $(\Li + \Rg) \cap \L^{\a,\i}$.
\end{itemize}

More details of the implementation of {METHOD 1} and {METHOD 2} can be
found in Appendix \ref{sec:appendix:grac}.

\subsubsection{Rank deficiency}
\label{sec:rank_def}
%
Let $I:=\#\L^i$ be the number of atoms in the interface and $R:=\#
\Rg$ the number of interacting sites. The number of unknowns
$C_{\ell;\rho,\vsig}$ is then $IR^2$. For Method 1, the number of
force balance equations is $\frac12 \# (\L^\i+\Rg) \times R$, and the
number of energy consistency equation is $2IR$. For method 2, we have
fewer force balance equations, while the number of constraints for
$C_{\ell,\rho,\vsig}$ is less than $\frac12 IR^2$. It is therefore
easy to see that the number of unknowns is much bigger than the number
of equations.

\subsubsection{Least squares computation of reconstruction parameters}
The references \cite{Shimokawa:2004,E:2006,PRE-ac.2dcorners} construct
various examples, where reconstruction parameters can be determined
analytically to satisfy the energy and force patch tests
\eqref{eq:E_pt_grac} and \eqref{eq:F_pt_grac}. Instead, we propose to
solve them numerically in a preprocessing step.

Comparing the number of equations against the number of free
parameters (see \S \ref{sec:rank_def}), we observe that, {\em if a
  solution to \eqref{eq:E_pt_grac} and \eqref{eq:F_pt_grac} exists,
  then it cannot be unique.} A natural idea, therefore, is to use a
least-squares approach,
\begin{equation}
  \label{eq:least_squares}
  \text{minimize } \sum_{\ell \in \L^\i} \sum_{\rho,\vsig \in
    \Rg(\ell)} |C_{\ell;\rho,\vsig}|^2 \quad \text{ subject to \eqref{eq:E_pt_grac} and \eqref{eq:F_pt_grac}.}
\end{equation}
We warn from the outset against using \eqref{eq:least_squares} and
explain in \S~\ref{sec:optim_coeffs} that error estimates for QNL type
a/c coupling schemes suggest a different selection principle.

Further, in \S~\ref{sec:stab}, we propose to add a stabilisation
mechanism to the interface site potentials that we previously explored
in \cite{2013-stab.ac}. Our subsequent numerical experiments in
\S~\ref{sec:numerical_tests} demonstrate that, in general, both of
these modifications are required to obtain satisfactory accuracy of
the a/c method.

\subsection{Consistency and Optimisation of $C_\ell$}
\label{sec:optim_coeffs}
In \cite[Thm. 6.1]{Or:2011a} it is shown that, under the assumptions
that $d = 2$ and that the atomistic region $\Omega^\a$ is connected
(and additional natural technical assumptions), any a/c coupling
scheme of the type \eqref{eq:generic_ac_energy} satisfying the force
and energy patch tests \eqref{eq:force_pt}, \eqref{eq:energy_pt}
satisfies a {\em first-order consistency estimate}: if $y = \yB$ in
$\L \setminus \Omega$ and if $\tilde{y}$ is an $H^2_{\rm
  loc}$-conforming interpolant of $y$, then
\begin{equation}
  \label{eq:consistency}
  \b\< \del E(y) - \del E^\ac(I_h y), u_h \b\> \leq C_1
  \| h \D^2 \tilde{y} \|_{L^2(\tilde{\Omega}^\c)},
\end{equation}
where $C_1$ is independent of $y$. (An improved result for a
specific variant of GRAC is also proven in \cite{PRE-ac.2dcorners}.) 

Of particular interest for the present work is the dependence on $C_1$
on the reconstruction parameters $C_\ell$, which we can obtain from
Equation (6.4) in \cite[Thm. 6.1]{Or:2011a} and a brief calculation:
\begin{align}
  \label{eq:C1_dependence}
  & C_1 \leq C_1' \, (1 + {\rm width}(\L^\i)) \, \sum_{\rho,\vsig \in
    \Rg} |\rho|\,|\vsig|\, M_{\rho,\vsig} + C_1''\\
  \notag
  \notag & \text{where} \quad M_{\rho,\vsig} = \max_{\ell \in \L^\i}
  \sum_{\tau,\tau' \in \Rg(\ell)} |V_{\tau,\tau'}(C_\ell \cdot Dy(\ell))| \,
  |C_{\ell;\tau,\rho}|\,|C_{\ell;\tau',\vsig}|.
\end{align}
$C_1'$ is a generic constant and $C_1''$ does not depend on the
reconstruction parameters.

The estimate \eqref{eq:C1_dependence} is of course an overestimation
that was convenient for the analysis, whereas intuitively one may
think of 
\begin{displaymath}
  M(\ell) := \sum_{\rho,\vsig} |\rho|\,|\vsig| \sum_{\tau,\tau'}
  \b|V_{\tau,\tau'}(C_\ell \cdot Dy(\ell))\b|  \,
  |C_{\ell;\tau,\rho}|\,|C_{\ell;\tau',\vsig}|
\end{displaymath}
to be a realistic ($\ell$-dependent) pre-factor. Suppose now that we
make the generic structural assumption (see App.B.2 in
\cite{OrtnerTheil2012}, where this is discussed for an EAM type
potential) that $|V_{\tau,\tau'}(C_\ell \cdot Dy(\ell))| \lesssim
\omega(|\tau|)\,\omega(|\tau'|)$, where $\omega$ has some decay that
is determined by the interaction potential, then we obtain that
\begin{align*}
  M(\ell) &\lesssim \sum_{\rho,\vsig} |\rho|\,|\vsig| \sum_{\tau,\tau'}
  \omega(|\tau|)\,\omega(|\tau'|) \,
  |C_{\ell;\tau,\rho}|\,|C_{\ell;\tau',\vsig}| \\
  &= \B(\sum_{\rho,\tau} |\rho| \omega(|\tau|)
  |C_{\ell;\tau,\rho}|\B)\,
  \B(\sum_{\vsig,\tau'} |\vsig| \omega(|\tau'|)
  |C_{\ell;\tau',\vsig}|\B) \\
  &= \B(\sum_{\rho,\tau} |\rho| \omega(|\tau|)
  |C_{\ell;\tau,\rho}|\B)^2.
\end{align*}
This indicates that, instead of $\|C\|_{\ell^2}$, we should minimise
$\max_{\ell \in \L^\i} \sum_{\rho,\tau} |\rho| \omega(|\tau|)
|C_{\ell;\tau,\rho}|$.  Since we do not in general know the generic
weights $\omega$, we simply drop them, and instead minimise
$\sum_{\rho,\tau} |C_{\ell;\tau,\rho}|$. Further, taking the maximum
of $\ell \in \L^\i$ leads to a difficult and computationally expensive
multi-objective optimisation problem. Instead, we propose to minimise
the $\ell^1$-norm of all the coefficients:
\begin{equation}
  \label{eq:ell1-min}
  \text{ minimise } \sum_{\ell \in \L^\i} \sum_{\rho,\vsig \in
    \Rg(\ell)} |C_{\ell;\rho,\vsig}| \quad \text{ subject to
    \eqref{eq:E_pt_grac} and \eqref{eq:F_pt_grac}.} 
\end{equation}

To justify the two rather significant simplifications, we observe
that, intuitively, the reconstruction coefficients at different sites
should take values of roughly the same order of magnitude. Further,
the weight factors coming from the interaction potential should not
play a big role since the reconstruction of each ``shell'' of
neighbours is in essence independent of the rest (due to the fact that
the reconstruction coefficients must also be valid for potentials with
smaller interaction neighbourhood). Finally, we remark that
$\ell^1$-minimisation tends to generate ``sparse'' reconstruction
parameters which may present some gain in computational cost in the
energy and force assembly routines for $E^\ac$.

\subsection{Stability and stabilisation}
\label{sec:stab}
In order to obtain an {\em energy norm error estimate} 
\begin{equation}
  \label{eq:err-est}
  \| D y_h - D y \|_{\ell^2} \leq C \B(\| h \D^2 \tilde{y}
  \|_{L^2(\tilde{\Omega}^\c)} + \epsilon^{\rm bc} \B),
\end{equation}
where $\epsilon^{\rm bc}$ is the error due to the artificial boundary
condition on $\partial\Omega$, we require a {\em best approximation
  error estimate}, the {\em consistency error estimate}
\eqref{eq:consistency}, and most crucially, a stability estimate of
the form
\begin{equation}
  \label{eq:stability}
  \< \ddel E^\ac(I_h y) u_h, u_h \> 
  \geq c_0 \| D u_h \|_{\ell^2}^2
\end{equation}
for some $c_0 > 0$, indepent of any approximation parameters. 

Estimates of the form \eqref{eq:stability} for any form of A/C
couplings in dimension greater than one are still poorly
understood. We refer to \cite{2013-stab.ac, 2012-MMS-bqcf.stab,
  LuMing:sharp} for some preliminary results. For our purposes, the
key observations from \cite{2013-stab.ac} are the following:
\begin{enumerate}
\item There exists no GRAC type a/c coupling for which
  \eqref{eq:stability} can be expected for general potentials $V$ and
  general boundary conditions $y_0$ {\em even} if $y$ itself is stable
  in the atomistic model.
\item By adding a stabilisation of the from $\kappa |D^2 y|^2$,
  with $\kappa$ sufficiently large, to the interface region,
  \eqref{eq:stability} can be expected. (We say ``expected'' instead
  of ``guaranteed'' since the proof of this statement in
  \cite{2013-stab.ac} is restricted to some specific interaction
  classes.)
\end{enumerate}

Thus, we shall consider also stabilised GRAC type couplings, where the
interface site potential is given by
\begin{equation}
  \label{eq:qnl_stab}
  \Phi_\ell^\i(y_h) := 
  V \b( C_\ell \cdot Dy_h(\ell) \b) + \kappa |D_{\rm nn}^2 y_h(\ell)|^2,
\end{equation}
where $\kappa \geq 0$ is a stabilisation parameter, and $|D_{\rm nn}^2
u_h(\ell)|^2$ is defined as follows: we choose $m \geq d$ linearly
independent ``nearest-neighbour'' directions $b_1, \dots, b_m$ in the
lattice, and denote
\begin{displaymath}
  \b|D_{\rm nn}^2 u_h(\ell)\b|^2 := \sum_{j = 1}^m \b| y_h(\ell+b_j) -
  2 y_h(\ell) + y_h(\ell-b_j) \b|^2.
\end{displaymath}
The reconstruction parameters $C_\ell$ are still determined according
to \eqref{eq:least_squares} or \eqref{eq:ell1-min}.

It is straightforward to see that the stabilisation does not generate
any ghost forces. That is, if the GRAC part of the potential,
$V(C_\ell \cdot Dy_h)$, satisfies the two patch tests
\eqref{eq:force_pt} and \eqref{eq:energy_pt}, then the stabilised
interface potential $\Phi_\ell^\i$ defined by \eqref{eq:qnl_stab} also
satisfies both patch tests.

\section{Numerical Tests}
\label{sec:numerical_tests}
\subsection{Model problems}
Our implementation is for the 2D triangular lattice $\mA \Z^2$ defined
by 
\begin{displaymath}
  \mA = \begin{pmatrix} 1 & \cos(\pi/3) \\ 0 & \sin(\pi/3) \end{pmatrix}.
\end{displaymath}
To generate a defect, we remove $k$ atoms
\begin{displaymath}
  \cases{ \L^{\defc}_{k} := \b\{ - (k/2+1) e_1, \dots , k/2 e_1\b\}, &
    \text{ if $k$ is even, } \\[1mm]
    \L^{\defc}_{k} :=  \b\{ - (k-1)/2 e_1, \dots, (k-1)/2 e_1 \b\}, & \text{ if $k$ is odd,}
  }
\end{displaymath}
to obtain $\L := \mA\Z^2 \setminus \L^{\defc}_k$. For small $k$, the
defect acts like a point defect, while for large $k$ it acts like a
small crack embedded in the crystal. In our experiments we shall
consider $k = 2, 11$.

We choose an elongated hexagonal domain $\Omega^\a$ containing $K$
layers of atoms surrounding the vacancy sites and the full
computational domain $\Omega$ to be an elongated hexagon containing
$N$ layers of atoms surrounding the vacancy sites; see Figure
\ref{fig:point_defects}(c) for an illustration. The domain
  parameters are chosen so that $N \approx K^2$. The finite element
  mesh is graded so that the mesh size function $h(x) = {\rm diam}(T)$
  for $T \in \T$ satisfies $h(x) \approx (|x|/K)^{3/2}$. These choices
  balance the coupling error at the interface, the finite element
  interpolation error and the far-field truncation error
  \cite[Sec.5.2]{2013-defects}. One then obtains
  \cite[Prop. 5.5]{2013-defects} under additional conditions on the
  stability of the method and the magnitude of the reconstruction
  parameters (we can verify both only {\it a posteriori}) that
\begin{equation}
  \label{eq:errest}
  \| \D y - \D y_h \|_{L^2} \leq C {\rm DOF}^{-1},
\end{equation}
where $y$ is identified with its P1 interpolant on the canonical
triangulation of $\L$ and ${\rm DOF}$ denotes the total number of
degrees of freedom (i.e. the number of atomistic sites $\L^{\a,\i}$
plus the number of finite element nodes).



The site energy is given by an EAM (toy-)model \eqref{eq:eam}, with 
\begin{align*}
  &\phi(r) = [e^{-2a(r-1)}-2e^{-a(r-1)}], \quad \psi(r) = e^{-br}, \\
  &F(\tilde{\rho})=c\b[(\tilde{\rho}-\tilde{\rho}_0)^2+(\tilde{\rho}-\tilde{\rho}_0)^4\b],
\end{align*}
with parameters $a = 4.4, b = 3, c=5, \tilde{\rho}_0 = 6e^{-b}$. The
interaction range is $\Nhd(\ell) = \L \cap B_2(\ell)$, i.e., next
nearest neighbors in hopping distance.

\subsubsection{Di-vacancy}
\label{sec:di-vacancy}
In the di-vacancy test two neighboring sites are removed, i.e., $k =
2$. We apply $3\%$ isotropic stretch and $3\%$ shear loading, by
setting
\begin{displaymath}
\mB := \begin{pmatrix} 1+s & \gamma_{\rm II} \\ 0 & 1+s \end{pmatrix}\cdot \mF_0.
\end{displaymath}
where $\mF_0 \propto I$ minimizes $W$, $s=\gamma_{\rm II}=0.03$.


\subsubsection{Micro-crack}
\label{sec:micro-crack}

In the microcrack experiment, we remove a longer segment of atoms,
$\L^{\defc}_{11}=\{-5e_1,\dots,5e_1\}$ from the computational
domain. The body is then loaded in mixed mode ${\rm I}$ \& ${\rm II}$,
by setting,
\begin{displaymath}
\mB := \begin{pmatrix} 1 & \gamma_{\rm II} \\ 0 & 1+\gamma_{\rm I} \end{pmatrix}\cdot \mF_0.
\end{displaymath}
where $\mF_0\propto I$ minimizes $W$, and $\gamma_{\rm I}=\gamma_{\rm
  II}=0.03$ ($3\%$ shear and $3\%$ tensile stretch).

\subsection{Methods}
We shall test the GRAC variants METHOD 1, METHOD 2 with both least
squares solution \eqref{eq:least_squares} and $\ell^1$-minimisation
\eqref{eq:ell1-min} to solve for the reconstruction parameters, and
with stabilisation parameters $\kappa = 0, 1$. The resulting methods
are denoted by {\tt M$i$-L$p$-S$\kappa$}, where $i \in \{1,2\}$, $p
\in \{2,1\}$, $\kappa \in \{0, 1\}$.

Some additional practical details for the implementation of METHOD 1
and METHOD 2 are described in Appendix \ref{sec:appendix:grac}.

We compare the GRAC methods with the five competitors previously
considered in \cite{2013-bqcfcomp, 2012-optbqce}:
\begin{itemize}
\item {\bf ATM}: full atomistic model is minimized with the constraint
  $y = \yB$ in $\L \setminus \Omega$; see also
  \cite[Sec. 4.1]{2013-defects}.
\item {\bf QCE}: original quasicontinuum method {\em without}
  ghost-force correction \cite{Ortiz:1995a}.
\item {\bf B-QCE, B-QCE+}: blended quasicontinuum method,
  implementation based on \cite{2012-optbqce}; B-QCE+ is a variant
  with highly optimised approximation parameters described in
  \cite[Sec. 4.3]{2013-bqcfcomp}.
\item {\bf QCF}: sharp-interface force-based a/c coupling
  \cite{Dobson:arXiv0903.0610}, formally equivalent to the
  quasi-continuum method {\em with} ghost-force correction
  \cite{Shenoy:1999a}.
\item {\bf B-QCF}: blended force-based a/c coupling, as described in
  \cite{2013-bqcfcomp}.
\end{itemize}

\subsection{Results}
Following \cite{2012-optbqce,2013-bqcfcomp} we present two
experiments, a di-vacancy ($k = 2$) and a ``micro-crack'' ($k =
11$). In the first experiment, we are able to clearly observe the
asymptotic behaviour of the a/c coupling schemes predicted in
\eqref{eq:errest}, while in the second experiment we observe a
significant pre-asymptotic regime where the prediction
  \eqref{eq:errest} becomes relevant only at fairly high DOF.

For both experiments we plot the absolute errors against the number of
degrees of freedom (DOF), which is propoertional to computational
cost, in the $H^1$-seminorm, the $W^{1,\infty}$-seminorm and in
the (relative) energy.

\begin{figure}
  \begin{center}
    \includegraphics[width=12.5cm]{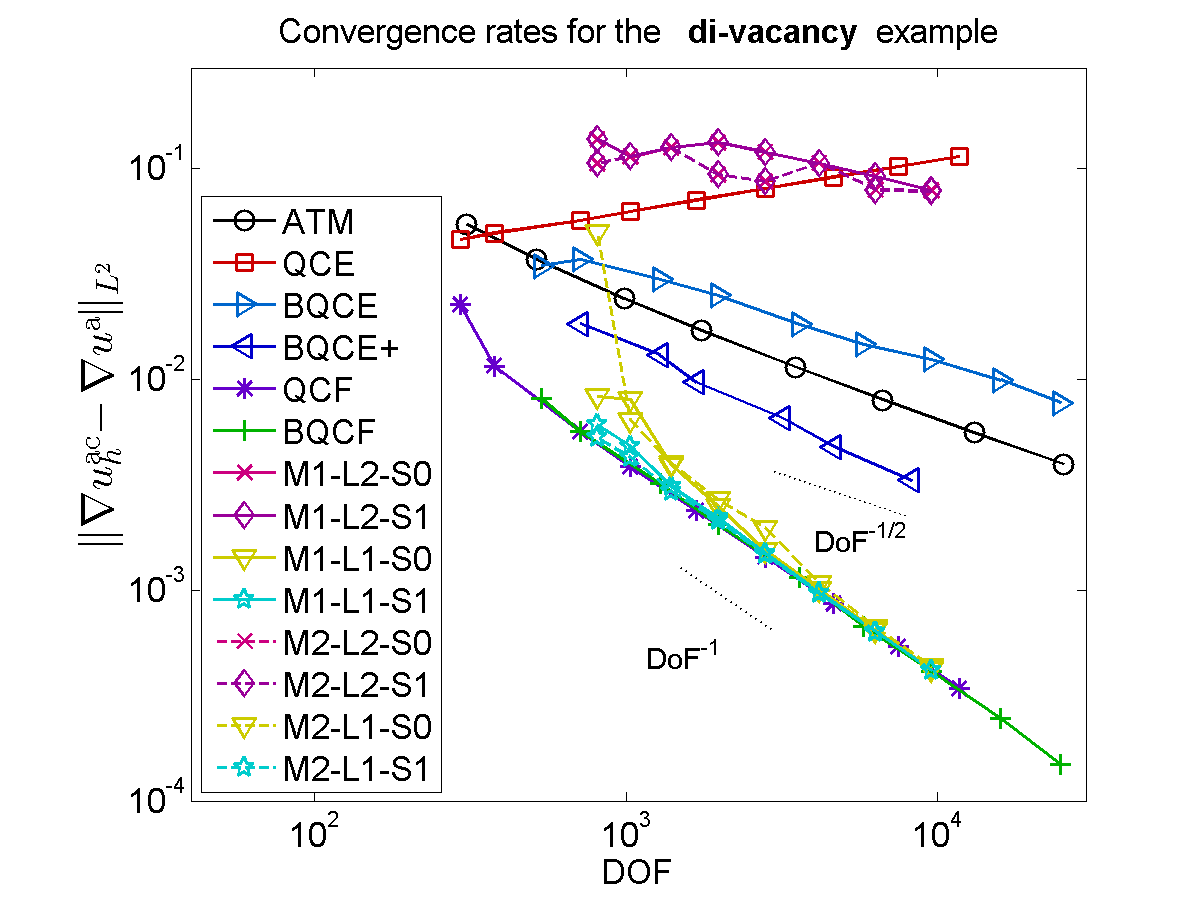}
    \caption{\label{fig:rate_divacancy_w12} Convergence rates in the energy-norm (the $H^1$-seminorm) for the divacancy benchmark problem described in Section \S~\ref{sec:di-vacancy} .}
  \end{center}
\end{figure}

\begin{figure}
  \begin{center}
    \includegraphics[width=12.5cm]{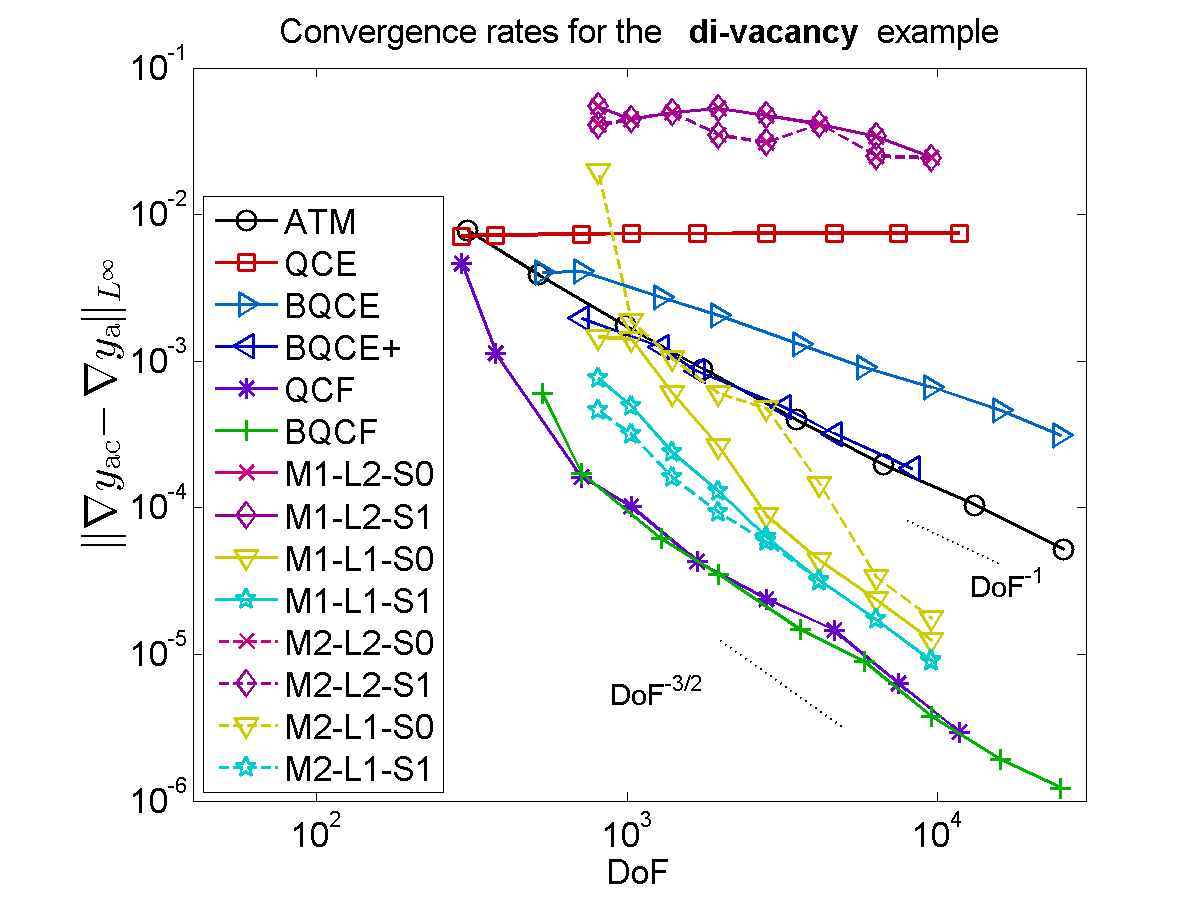}
    \caption{\label{fig:rate_divacancy_w1inf} Convergence rates in the $W^{1,\infty}$-seminorm for the divacancy benchmark problem described in Section \S~\ref{sec:di-vacancy} .}
  \end{center}
\end{figure}

\begin{figure}
  \begin{center}
    \includegraphics[width=12cm]{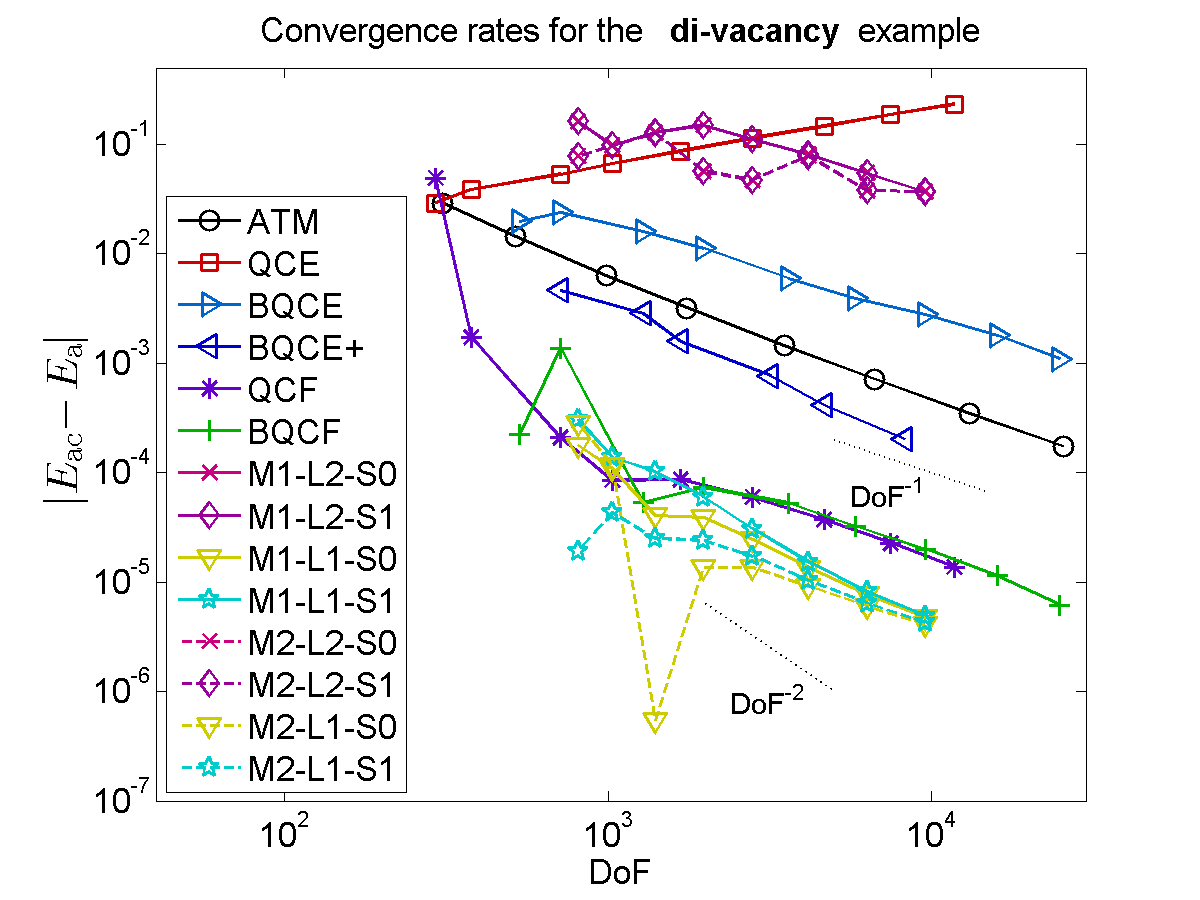}
    \caption{\label{fig:rate_divacancy_E} Convergence rates in the
      {relative energy} for the divacancy benchmark problem described
      in Section \S~\ref{sec:di-vacancy} .}
  \end{center}
\end{figure}

\begin{figure}
  \begin{center}
    \includegraphics[width=12cm]{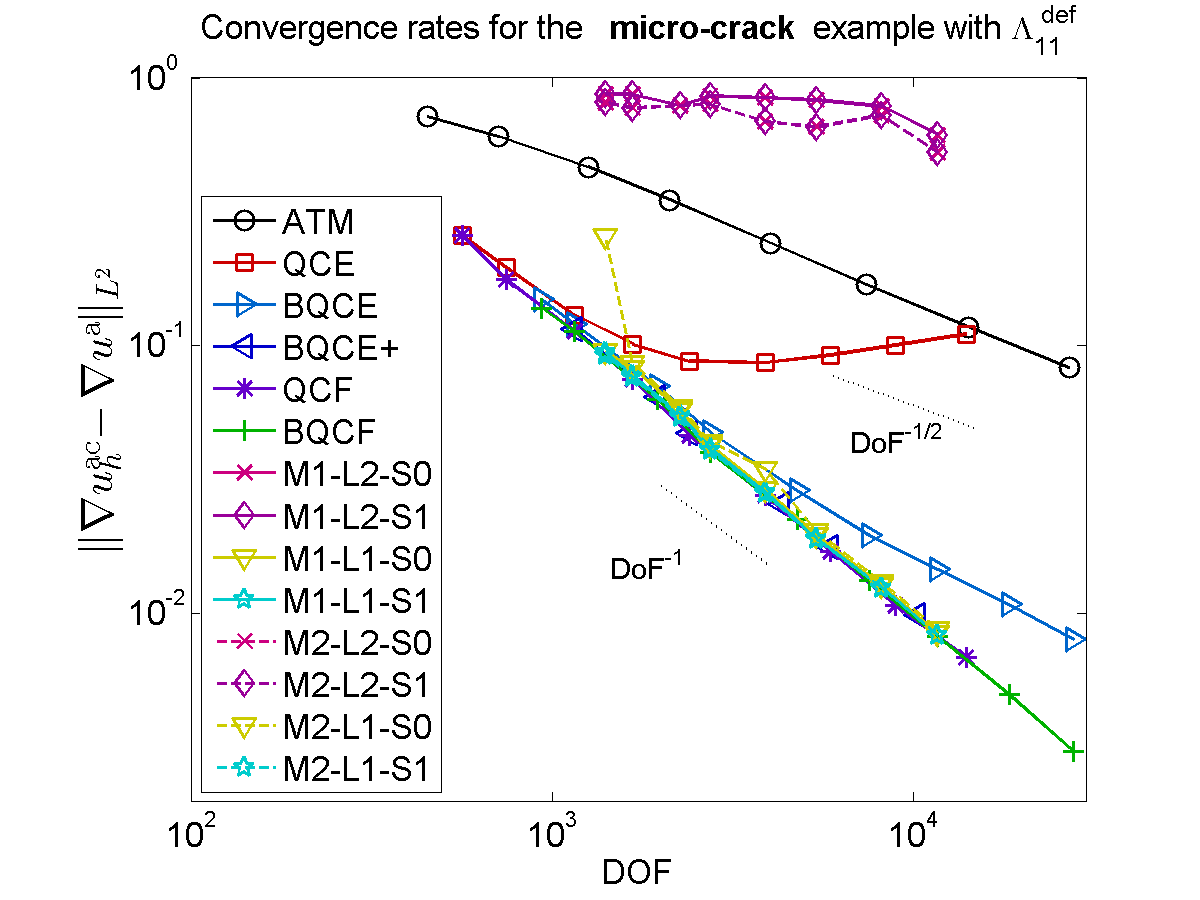}
    \caption{\label{fig:rate_microcrack11_w12} Convergence rates in the energy-norm (the $H^1$-seminorm) for the microcrack benchmark problem with $\L^{\defc}_{11}$ described in Section \S~\ref{sec:micro-crack} .}
  \end{center}
\end{figure}

\begin{figure}
  \begin{center}
    \includegraphics[width=12cm]{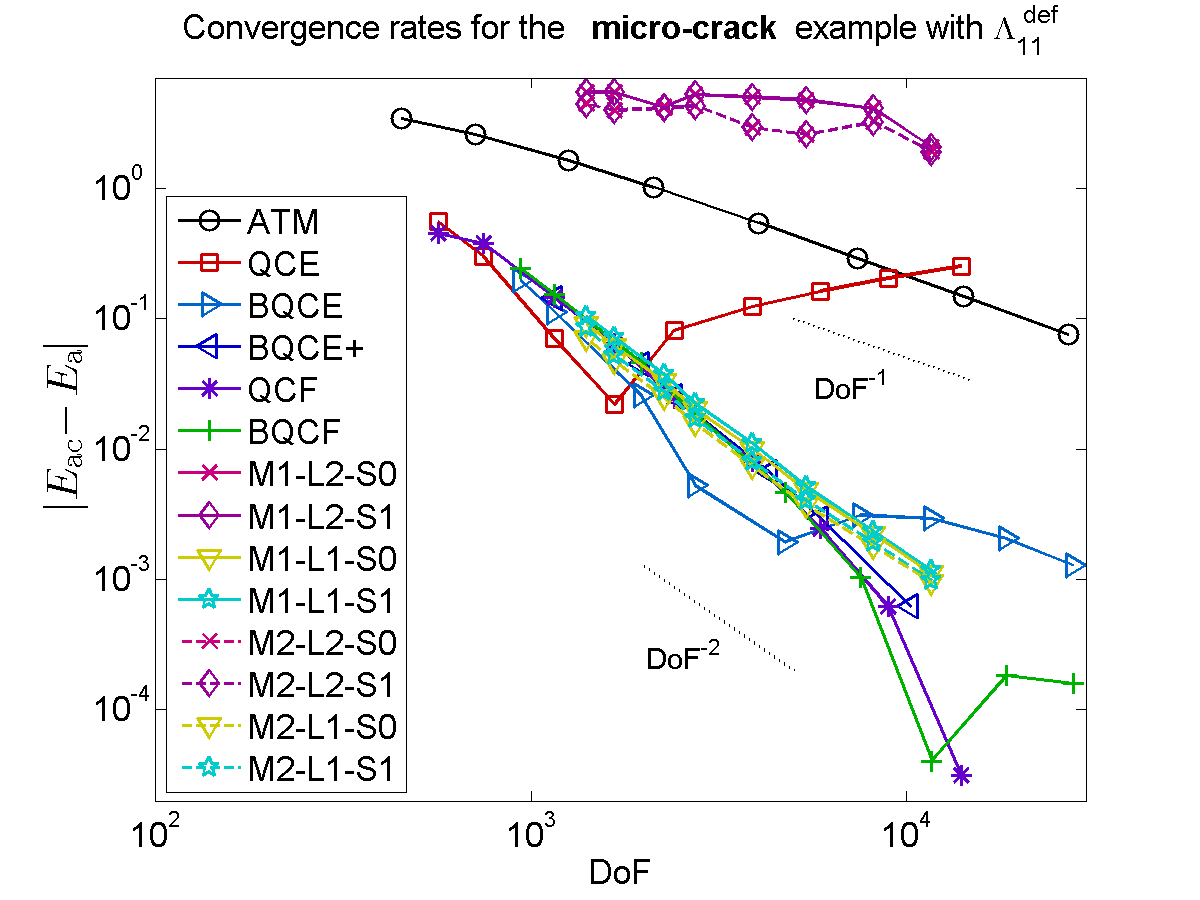}
    \caption{\label{fig:rate_microcrack11_E} Convergence rates in the
      {relative energy} for the microcrack benchmark problem with
      $\L^{\defc}_{11}$ described in Section \S~\ref{sec:micro-crack}
      .}
  \end{center}
\end{figure}

The results are shown in Figures \ref{fig:rate_divacancy_w12},
\ref{fig:rate_divacancy_w1inf} and \ref{fig:rate_divacancy_E} for the
divacancy problem and in Figures \ref{fig:rate_microcrack11_w12},
\ref{fig:rate_microcrack11_w1inf} and \ref{fig:rate_microcrack11_E}
for the micro-crack problem.




\subsubsection{Effect of $\ell^1$-minimisation}
In all error graphs we observe that computing the reconstruction
coefficients via least-squares ($\ell^2$-minimisation) leads to large
errors in the computed solution and likely even lack of
convergence. Stabilisation does not remedy this, which indicates that
the issue indeed lies in the consistency error. By contrast, using
\eqref{eq:ell1-min} ($\ell^1$-minimisation) to compute the
reconstruction parameters leads to errors that are competitive with
the provably quasi-optimal schemes QCF and B-QCF.

\subsubsection{Effect of  stabilisation}
\begin{figure}
  \begin{center}
    \includegraphics[width=12cm]{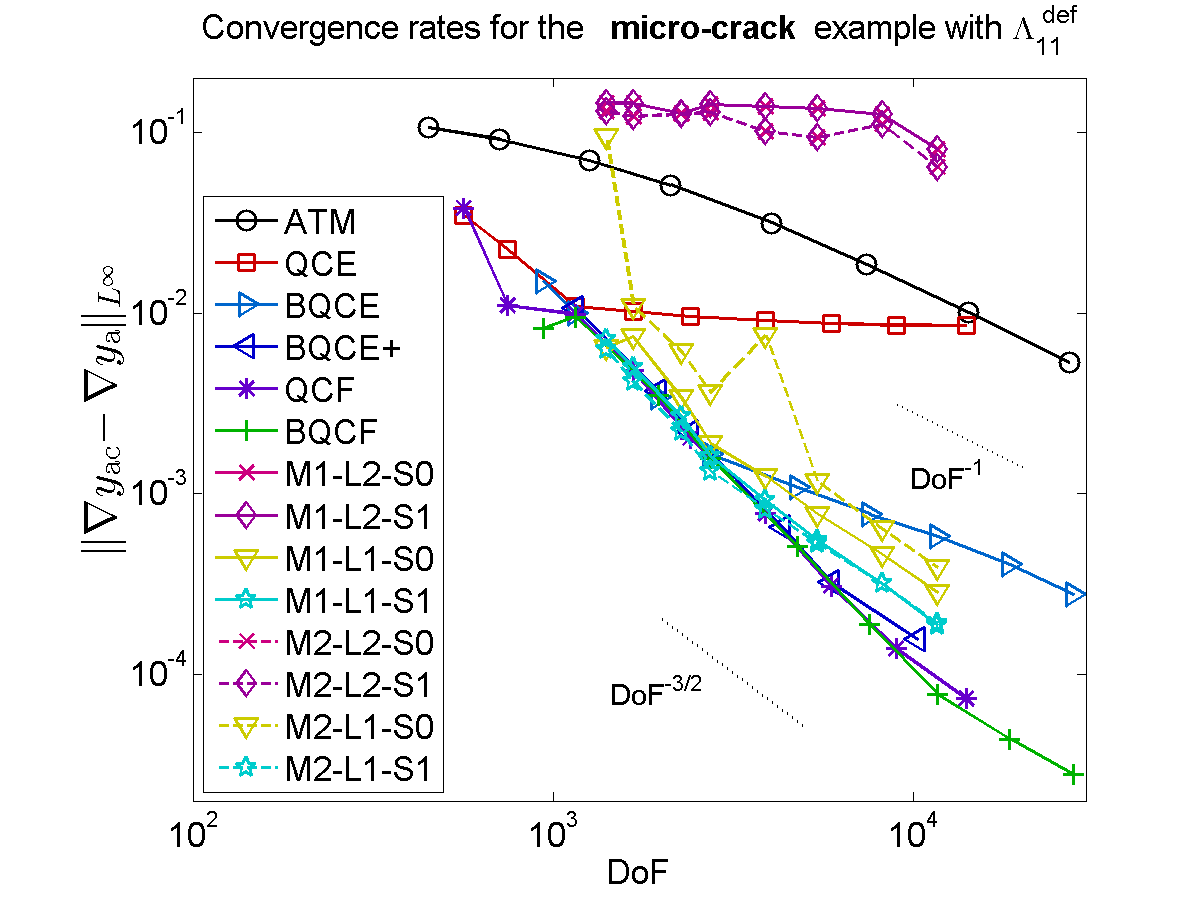}
    \caption{\label{fig:rate_microcrack11_w1inf} Convergence rates in
      the {$W^{1,\infty}$-seminorm} for the microcrack benchmark
      problem with $\L^{\defc}_{11}$ described in Section
      \S~\ref{sec:micro-crack} .}
  \end{center}
\end{figure}

If no stabilisation is used ($\kappa = 0$), then all error graphs
display large errors in a pre-asymptotic regime and in some cases,
most pronounced in Figure \ref{fig:rate_microcrack11_w1inf},
non-monotone convergence history.

Adding the stabilisation by setting $\kappa = 1$ the $H^1$ and
$W^{1,\infty}$ errors are reduced in both examples, indeed
significantly so in the important pre-asymptotic regime, and the
oscillations in the convergence history are removed. With
stabilisation the convergence rates predicted in
\cite[Sec. 5.2]{2013-defects} are clearly observed.




\subsubsection{Comparison of a/c couplings}
%
%
%
%
In all error graphs we clearly observe the optimal convergence rate of
GRAC ({\tt M$i$-L1-S1} variants) among the tested energy-based methods
(ATM, QCE, B-QCE, B-QCE+, GRAC). Indeed, the errors are even
competitive with the quasi-optimal force-based schemes (QCF, B-QCF):
for $H^1$ errors they are essentially comparable, for
$W^{1,\infty}$ errors the force-based schemes are only better by a
moderate constant factor, while for the energy errors the GRAC methods
are optimal. (Note that, for QCF we evaluate the QCE energy and for
B-QCF we evaluate the B-QCE energy.)


\clearpage

\section{Conclusion}
\label{sec:discuss}
We have succeeded in presenting the first patch test consistent
energy-based atomistic-to-continuum coupling formulation, GRAC, which
is applicable to general a/c interface geometries and general
(short-ranged) many-body interactions, and demonstrated its potential
in a 2D implementation.

We have discussed the critical issues of $\ell^1$-minimisation and of
stabilisation, and have demonstrated that our final formulations yield
an energy-based a/c coupling that is optimal among the energy-based
methods we tested, which represent a fairly generic sample, and are
even competitive compared against the quasi-optimal force-based
coupling schemes.

While the construction of the GRAC scheme is involved, it has the
  advantage that no additional approximation parameters (e.g., the
  blending function $\beta$ in the B-QCE and B-QCF schemes
  \cite{2012-optbqce,2013-bqcfcomp}) must be adapted to the problem at
  hand.

The main challenge that requires additional work is the complexity of
the precomputation of the reconstruction parameters, which may become
prohibitive for wider interaction stencils, in particular in 3D. It
may then become necessary to make further simplifications such as the
ones we made in METHOD 2, in order to substantially reduce the
computational cost and storage to compute these parameters.

From a theoretical perspective the main open problem is to prove that
the geometric consistency equations \eqref{eq:E_pt_grac} and
\eqref{eq:F_pt_grac} always have at least one solution. We can, at
present, provide no analytical evidence to support this claim,
however, we have so far not encountered a situation where a solution
could not be computed numerically.

Finally, we remark that the consistency of the GRAC scheme is
  still not entirely settled. First-order consistency is only proven
  in 1D and in 2D under the restrictive assumption that the atomistic
  region is connected \cite{Or:2011a, LZ:2014}.

%
%

\clearpage 

\section{Appendix}
\label{sec:appendix}

\subsection{First variation of $\Eac$}
\label{sec:app_delEac}
The following calculations provide the details for the computation of
$\del\Eac$ in \S~\ref{sec:force_pt_lineqn}.

\subsubsection{Atomistic component}
\begin{align*}
  \<\del\Ea, u\>  = & \sum_{\ell\in\La}\sum_{\rho\in\Rg} \D_\rho V D_\rho u(\ell)\\
  = & \sum_{\ell\in\La}\sum_{\rho\in\Rgp} \D_\rho
  V(u(\ell+\rho)-u(\ell)) +
  \D_{-\rho}V(u(\ell-\rho)-u(\ell))\\
  = &\sum_{\ell\in\La}\sum_{\rho\in\Rgp}\D_\rho V(u(\ell+\rho)-u(\ell-\rho))\\
  = &\sum_{\ell+\rho\in\La,\rho\in\Rgp}\b[\D_\rho V\cdot u(\ell)\b] - \sum_{\ell-\rho\in\La, \rho\in\Rgp}\b[\D_\rho V \cdot u(\ell)\b]\\
\end{align*}

\subsubsection{Interface component}
\begin{align*}
	\<\del\Ei, u\>  = &\sum_{\ell\in\Li}\omega_\ell^\i\<\del V\b((\sum_{\vsig\in\Rg} C_{\ell;\rho,\vsig} D_\vsig y(\ell))_{\rho\in\Rg}\b), u\>\\
								  = &\sum_{\ell\in\Li}\omega_\ell^\i\sum_{\rho\in\Rg}\sum_{\vsig\in\Rg}C_{\ell;\rho,\vsig} \D_\rho V D_\vsig u(\ell)\\
								  = & \sum_{\ell\in\Li}\omega_\ell^\i\sum_{\rho\in\Rgp}\sum_{\vsig\in\Rg}(C_{\ell;\rho,\vsig}-C_{\ell;-\rho,\vsig})\D_\rho V (u(\ell+\vsig)-u(\ell))\\
								  = &\sum_{\ell-\vsig\in\Li, \vsig\in \Rg}\omega_{\ell-\vsig}^\i\sum_{\rho\in\Rgp}(C_{\ell-\vsig;\rho,\vsig}-C_{\ell-\vsig;-\rho,\vsig})\b[\D_\rho V \cdot u(\ell)\b]\\
				& -\sum_{\ell\in\Li}\omega_\ell^\i\sum_{\rho\in\Rgp}\sum_{\vsig\in\Rg}(C_{\ell;\rho,\vsig}-C_{\ell;-\rho,\vsig})\b[\D_\rho V \cdot u(\ell)\b]
\end{align*}

\subsubsection{Cauchy--Born component}
\begin{align*}
	\<\del\Ec, u\>  = &\sum_T v_T\<\del W, u\>\\
								  = &\sum_T \frac{v_T}{|\vor|}\<\del V\b((\D_Ty\cdot \rho)_{\rho\in\Rg}\b),u\>\\
								  = &\sum_T \frac{v_T}{|\vor|}\sum_{\rho\in\Rg}\D_\rho V\D_T u\cdot \rho\\
								  = &\sum_T \frac{v_T}{|\vor|}\sum_{\rho\in\Rg}\D_\rho V\sum_{i=1}^3u^T_i\D_T\phi_i^T\cdot \rho\\
								  = &\sum_T \frac{v_T}{|\vor|}\sum_{\rho\in\Rgp}2\D_\rho V\sum_{i=1}^3u^T_i\D_T\phi_i^T\cdot \rho\\
								  = &\sum_T\sum_{\rho\in\Rgp}\sum_{i=1}^3 2\frac{v_T}{|\vor|}\D_T\phi_i^T\cdot \rho \b[\D_\rho V \cdot u^T_i\b]
\end{align*}

\subsection{Setup of the geometric consistency equations}
\label{sec:appendix:grac}
We now introduce additional details for implementing the GRAC
formulation in \eqref{eq:generic_ac_energy}. This gives
further concrete details on how to setup the geometric consistency
equations \eqref{eq:E_pt_grac} and \eqref{eq:F_pt_grac} specifically
for the triangular lattice. The process that we propose is, however,
more generally applicable. Here, the interface region is $r$ layers of 
atoms around $\L^\a$, and $r$ is the radius of interaction range $\Rg$ 
in terms of hopping distance. We describe the process only for METHOD 1,
as the one for METHOD 2 is very similar.

To satisfy the force patch test consistency equation, in the nearest
neighbor case we considered in \cite{PRE-ac.2dcorners} we take the
following strategy, where the reconstruction parameters $C_\ell$ are
extended to $\L \setminus \L^\i$ by
  \begin{displaymath}
    C_\ell = \cases{C^\a, & \ell \in \L^\a, \\
      C^\c, & \ell \in \L \setminus \L^{\a,\i},}
  \end{displaymath}
Define
the six nearest-neighbour lattice directions by $a_1 := (1, 0)$, and
$a_j := \mQ_6^{j-1} a_1$, $j \in \Z$, where $\mQ_6$ denotes the
rotation through angle $2\pi/6$ and we note that $a_{j+3} = -
a_j$.  Then $C^\a$ is given by
  $C^\a_{i,j} = \delta_{i,j}$, $C^\c$ is given by $C^\c_{i,j} =
  \frac23\delta_{i,j}+\frac13\delta_{i,j+1}+\frac13\delta_{i,j-1}$, $i,j=1,\dots,6$,
  where $\delta$ is the Kronecker delta function.

  The argument employed in \cite[Lemma 3.2]{PRE-ac.2dcorners} can be
  extended to longer range interactions. There exist matrices
  $C_\ell^\c$ such that, upon defining $\Psi_\ell(y) := V( C_\ell^\c
  \cdot Dy(\ell))$, we have
\begin{equation}
\label{eq:cbcoefs}
  \< \del \Psi_\ell(\mF x), v \> = \int_{{\rm vor}(\ell)} \partial
  W(\mF) : \D v(x) \dx \qquad \forall \ell \in \L \cap \Omega^\c,
\end{equation}
that is, under uniform deformation, the forces generated by the
Cauchy--Born site potential $\int_{{\rm vor}(\ell)} W(\D y) \dx$ are
the same as those of $\Phi_\ell$. 

Carrying this out in practise requires that several layers of atoms
surrounding $\L^{\a,\i}$, denoted by $\L^\c$, coincide with the finite
element nodes in that region.
%
Upon choosing $\Tmu$ to be a uniform partition over $\L^\c$, these
parameters can be computed analytically.  The details are shown in
Appendix \ref{sec:appendix:cc} for next nearest neighbor interactions.

Upon defining the coefficients for the atomistic and continuum region,
we can use Proposition \ref{prop:coefficients} to compute unknown
parameters $C_{\ell; \rho,\vsig}$.




\subsection{Determination of the coefficients $C^\c$ for next nearest neighbor interaction.}
\label{sec:appendix:cc}
We now calculate the coefficients $C^\c$ from the equation
\eqref{eq:cbcoefs}. On the canonical triangular mesh induced by $\L$, let
\begin{displaymath}
  \displaystyle V^\c_\ell = \frac{1}{6}\sum_{T \owns \ell} V(D_T u)
\end{displaymath}
be the Cauchy-Born site energy with respect to $\ell\in\L$. As 
the six nearest-neighbour lattice directions are defined in Section \ref{sec:appendix:grac}. The second nearest-neighbour lattice directions can be expressed as $a_{2j+5} = 2a_j$, $a_{2j+6} = a_j + a_{\shift(j)}$, $j=1,\dots,6$, , where $\shift\{1,2,3,4,5,6\}=\{2,3,4,5,6,1\}$. Therefore $a_j$'s, $j=1,\dots,18$ form the the interaction
range $\Rg$ for next nearest neighbor interactions.

$V_\ell^\c$ only depends on the first 6 variables $D_iy$ of $V$, a direct calculation shows that 
\begin{displaymath}
	\pp_1 V_\ell^\c = \frac13\pp_2 V + \frac23\pp_1V + \frac13\pp_6 V + \frac23\pp_9 V + \frac23 \pp_9 V + \pp_8 V + \frac{4}{3}\pp_7 V + \pp_{18}V + \frac23\pp_{17}V
\end{displaymath}
and similiarly for $\pp_i V_\ell^\c$ with $i = 2,\dots, 6$. 

Now we can write down the modified potential $\Psi_\ell$ defined in \eqref{eq:cbcoefs}, which generate the same force for arbitrary uniform deformations. In the following expression of $\Psi_\ell$, for $i = 1, \dots, 6$, $D_i y$  are abbreviated by $D_i$,
\begin{align*}
	\Psi_\ell = &V\b(\frac23 D_1+\frac13 D_2+\frac13 D_6,\frac23 D_2+\frac13 D_1+\frac13 D_3,\frac23 D_3+\frac13 D_2+\frac13 D_4,\\
								 &\frac23 D_4+\frac13 D_3+\frac13 D_5,\frac23 D_5+\frac13 D_4+\frac13 D_6,\frac23 D_6+\frac13 D_5+\frac13 D_1,\\
								 &\frac{4}{3}D_1+\frac23 D_2+\frac23 D_6, D_1+D_2, \frac{4}{3}D_2+\frac23 D_1+\frac23 D_3, \\
								 &D_2+D_3,\frac{4}{3}D_3+\frac23 D_2+\frac23 D_4,D_3+D_4\\
								 &\frac{4}{3}D_4+\frac23 D_3+\frac23 D_5, D_4+D_5, \frac{4}{3}D_5+\frac23 D_4+\frac23 D_6, \\
								 &D_5+D_6,\frac{4}{3}D_6+\frac23 D_5+\frac23 D_1,D_6+D_1\b)
\end{align*}
Hence the coefficients $C_\ell^\c$ can be drawn from above expression by using $\Psi_\ell = V( C_\ell^\c \cdot Dy(\ell))$.

\bibliographystyle{plain}
\bibliography{qc}
\end{document}